\documentclass{article}
\usepackage{graphicx} 
\usepackage{amsmath, amssymb, amsthm}
\usepackage{graphicx}
\usepackage{bm}
\usepackage{upgreek}
\usepackage{enumitem} 
\usepackage{hyperref}
\hypersetup{colorlinks,linkcolor={blue},citecolor={red},urlcolor={blue}}
\usepackage{cleveref}
\usepackage{mathtools}
\usepackage{float}
\usepackage[a4paper, total={6in, 10in}]{geometry}
\usepackage{enumitem} 
\usepackage{authblk}

\newtheorem{thm}{Theorem}[section]
\newtheorem{cor}{Corollary}[section]
\newtheorem{lemma}{Lemma}[section]

\newtheorem{Aassumption}{A-Assumption}
\newtheorem{Bassumption}{B-Assumption}

\newtheorem{prop}{Proposition}

\newtheorem{assmp}{Assumption}

\newtheorem{note}{Note}

\title{Ignoring Blow-Ups for Approximations}
\author{Tim Johnston\footnote{ESSEC Paris, 	
timothy.johnston@essec.edu}}
\date{October 2026}

\begin{document}
\maketitle
\section{Introduction}
\begin{abstract}
It is well known that for SDEs with non-Lipschitz coefficients, the Euler scheme can fail to converge in an $L^p$ sense. It is also widely observed that convergence fails only due to rare explosive trajectories, that distort the expectation. In this work we demonstrate that almost any method of removing these rare explosive trajectories leads to convergence in a suitable strong and weak sense. We apply this framework to coefficients that grow faster than linearly at infinity and Lennard-Jones type coefficients that explode like a reciprocal power on a compact set.
\end{abstract}
In this article we consider the simulation of SDEs taking values in $\mathbb{R}^d$ of the form
\begin{equation}\label{eq: main equation}
    dX_t=b(X_t)dt+dW_t,
\end{equation}
where the initial condition $X_0$ is a $\mathcal{F}_0$-measurable random variable independent of the driving noise, and conditions on the coefficients will be made more explicit further on. Whilst it is well known that explicit Euler-type approximations work well for Lipschitz coefficients \cite{kloeden2013numerical}, it has additionally been observed that if the coefficients grow superlinearly, see for example dynamics like \eqref{eq: example SDE superlinear}, then rare explosive trajectories can lead to non-convergence (or divergence) in $L^p$ as the step-size tends to $0$, see \cite{articlediv, 10.1214/12-AAP890}. This lead to the development of tamed or truncated schemes, which get around this issue by reducing large values of the drift coefficient as a function of step-size \cite{da2b46b6-84fa-3eed-be6d-20c143c60df1, articletamed, MAO2015370}. This development has been fruitful, and analougues of this idea have been extended to SPDEs \cite{0329d33d261740b59da99f06a80a9570, Brehier2022}, other types of noise \cite{doi:10.1137/151004872, Qi2026}, higher order schemes \cite{ChamanKumar2017, Wang2013}, sampling algorithms \cite{BROSSE20193638, theopouladiscont, doi:10.1137/22M1514283}, and interacting particle systems \cite{Johnston2025, Kumar2020WellposednessAT, soni2025tamedeulerapproximationfully}. An alternative method using a pathwise variant of the unadjusted Barker proposal algorithm \cite{10.1093/imanum/draf147} was recently proposed in \cite{iguchi2026pathwiseskewsymmetricdiscretisationsdes}. Furthermore, certain other recent works extend these ideas to coefficients that may blow up to infinity on a compact set \cite{Johnston2024ATE, 10.1214/24-AOP1750}, reflecting models commonly used in molecular dynamics \cite{29acd3d494044594aea0829ef236aad6}. 

It should be noted, however, that blow-ups (sometimes tens or hundreds of order of magnitude larger than the remainder of the values) are extremely rare when the step-size is sufficiently small. This is demonstrated for instance in the simulation section of \cite{articlediv}. One may wonder in this case what the harm is in simply ignoring these rare explosions, which appear as statistical outliers. 

In this article, we show that the size of the set of explosive trajectories for the Euler scheme is exponentially small in the reciprocal of the step-size. As a corollary to this, we demonstrate that almost any method for removing (polynomially) `large' trajectories restores strong convergence, and conditioning on this by rejecting such trajectories ensures weak convergence (in Wasserstein distance).

Similar ideas were explored in \cite{5aa5837d-f9a1-33da-b667-726b0d620667}, where all trajectories leaving a certain ball are discarded, and weak convergence is shown up to an error. Also related is the result of \cite{bourabee2010} on the Metropolis-Hastings adjusted variant of the Euler scheme (MALA), which similarly shows that super-linearities only cause issues up to exponentially small terms. In our work we show the error from \cite{5aa5837d-f9a1-33da-b667-726b0d620667} can in fact be removed, and that we can extend the results of \cite{5aa5837d-f9a1-33da-b667-726b0d620667} and \cite{bourabee2010} to SDEs with singular coefficients that take values in an open set $D\subset \mathbb{R}^d$. 

Our Lyapunov function approach to the regularity of the coefficients has aspects in common with \cite{5aa5837d-f9a1-33da-b667-726b0d620667}, as well as \cite{70a33031-e55f-31e6-80ab-f131fb94c87a}. However, we use a Lyapunov function to control the \textit{regularity} of the drift coefficient rather than the \textit{norm}. Therefore, by controlling this Lyapunov function for the scheme and true solution, we demonstrate that the dynamics `mostly' stay in a region where $b$ is `almost' Lipschitz. This solidifies the intuition that many classes of non-Lipschitz coefficients are `almost Lipschitz', in the sense that after explosive trajectories are removed both theory and practise is similar to as in the Lipschitz case. This is distinct from other classes of non-Lipschitz coefficients whose difficulty cannot be so simply overcome, studied for instance in \cite{arnulf2015, gerencser2017, MULLERGRONBACH20181013}, as well as discontinuous coefficients, for which despite being amenable to Euler-type methods \cite{Leobacher2016multi, articlemg, 10.1214/20-EJP479}, behave in other ways quite dissimilarly to Lipschitz coefficients in the context of Euler approximations \cite{10.1214/22-AAP1837}.

\section{Removing Trajectories}

In this paper we consider the Euler discretisation of \eqref{eq: main equation}, that is for $n\in \mathbb{N}$
\begin{equation}\label{eq: Euler scheme}
    dX^n_t=b(X^n_{\kappa_n(t)})dt+dW_t,
\end{equation}
where $\kappa_n(t)=\frac{\lfloor nt\rfloor}{n}$ is the projection backwards onto the grid $\{0,1/n,2/n,...\}$. One notes that the step-size of the Euler-type approximation $X^n_t$ is therefore $1/n$. In cases where the regular Euler scheme \eqref{eq: Euler scheme} is not stable as an estimator of \eqref{eq: main equation}, the tamed approximation paradigm involves finding new step-size dependent coefficients $b_n:\mathbb{R}^d\to\mathbb{R}^d$ and considering the approximation
\begin{equation}\label{eq: tamed Euler scheme}
    dX^n_t=b_n(X^n_{\kappa_n(t)})dt+dW_t.
\end{equation}
Here it is necessary that $b_n\to b$ in some suitable sense as $n\to\infty$, and $b_n$ can be given for instance like $b_n:=\frac{b}{1+n^{-1/2}b}$. Our aim shall be to show that, for a wide class of SDEs, one can recover convergence by simply \textit{ignoring} rare explosive trajectories. In particular for \textit{any} choice of $C,h>0$, one may recover strong convergence at rate $1$ by replacing with $0$ any trajectory that exceeds $Cn^r$ in Euclidean norm at grid-points. That is, for every $T>0$, in order to approximate the true solution $(X_t)_{t\in [0,T]}$ one may choose
\begin{equation}\label{eq: truncated approximation}
    Y^n_t=X^n_t\cdot 1_{\{\sup_{t\in [0,T]}\lvert X^n_{\kappa_n(t)}\rvert \geq Cn^h\}},
\end{equation}
which we demonstrate converges in $L^2$. This is the content of Corollary \ref{cor: strong convergence}. Furthermore, by conditioning on the set $\{\sup_{t\in [0,T]}\lvert X^n_{\kappa_n(t)}\rvert \leq Cn^r\}\subset \Omega$ we may recover weak convergence in the Wasserstein distance
\begin{equation}\label{eq: wasserstein}
    \mathcal{W}_2(\mu,\nu):=\inf_{(X,Y)\in T(\mu,\nu)}\sqrt{E\lvert X-Y\rvert^2},
\end{equation}
where $\mu,\nu$ are measures on a polish space $\mathcal{X}$ and $T(\mu,\nu)$ is the set of couplings on $\mathcal{X}\times \mathcal{X}$ for which $X\sim \mu$ and $Y\sim \nu$. This is the content of Corollary \ref{cor: weak convergence}.

These results are both corollaries of our main result which says that explosive trajectories are \textit{exponentially rare} in $n$. Therefore, any method which replaces these trajectories with something integrable, as well as which doesn't significantly change most other trajectories, should recover convergence.

In this paper we track \textbf{only dependence on the step-size parameter $n\in \mathbb{N}$}
\begin{note}
    Let $c>0$ denote a generic constant changing from line to line that depends on every global parameter (i.e. $b$, $d$, $D$,  all the constants listed in the assumptions, the terminal time $T>0$) \textbf{except} the stepsize parameter $n\in \mathbb{N}$.
\end{note}

\subsection{Assumptions}\label{subsec: assumptions}
We consider the equation \eqref{eq: main equation} which takes values in $\mathbb{R}^d$, driven by a $d$-dimensional Brownian motion defined on a filtered probability space $(\Omega, P, (\mathcal{F}_t)_{t\geq 0})$. We assume the solution takes values in an open set $D\subset\mathbb{R}^d$, the drift coefficient $b:D\to\mathbb{R}^d$ is a function defined on this open set (this shall be enforced by our assumptions, see Lemma \ref{lemma: exponential moments for true solution}). In order to ensure the numerical scheme is well-defined everywhere, it shall be convenient to extend this to a function $b:\mathbb{R}^d\to\mathbb{R}^d$ by letting $b$ be $0$ everywhere outside of $D\subset\mathbb{R}^d$. Since we shall ignore every trajectory of the scheme which leaves $D$, this distinction shall not be important. Therefore by an abuse of notation we refer to both functions as $b$, so there are no issues with writing \eqref{eq: Euler scheme}.

Key to our framework shall be a Lyapunov function $V:D\to [1,\infty)$ that we use to control the regularity of our drift coefficient $b$. We assume of $b, V, D$ and $X_0$ the following

\begin{assmp}[Quasi-Lipschitzness]\label{assmp1: Lipschitz like behaviour}
    The set $D\subset \mathbb{R}^d$ is convex and for all $x,y \in D$ one has
    \begin{equation}
        \lvert b(x)-b(y)\rvert\leq c(V(x)+V(y))\lvert x-y\rvert.
    \end{equation}
\end{assmp}

\begin{assmp}[Quasi-Smoothness]\label{assmp2: higher smoothness}
    One has that $b\in C^1(D,\mathbb{R}^d)$, and there exists $q>0$ such that for all $x,y \in D$ one has
    \begin{equation}
        \lvert b(x)-b(y)-\nabla b(x)(x-y) \rvert\leq c(V(x)^q+V(y)^q)\lvert x-y\rvert^2.
    \end{equation}
\end{assmp}
\begin{assmp}[Integrability of Moments]\label{assmp3: moment assmp}
    Suppose that $V\in C^2(D,[1,\infty))$ is twice differentiable, and in the case $D\not=\mathbb{R}^d$ one has $V(y)\to\infty$ as $D\ni y\to x\in \partial D$. Additionally, assume there exist $\alpha, \beta>0$ such that for every $r>0$ there exists $c_r>0$ such that
    \begin{equation}
  \mathcal{L}e^{\beta V^\alpha} +re^{\beta V^\alpha}\leq c_r,
    \end{equation}
   where $\mathcal{L}K:=\langle \nabla K , b\rangle +\frac{1}{2}\Delta K$ is the generator of \eqref{eq: main equation}. Furthermore, the initial condition $X_0$ is $\mathcal{F}_0$-measurable and satisfies
      \begin{equation}
E[e^{\beta V(X_0)^\alpha}]\leq c.
    \end{equation}
\end{assmp}
    \begin{assmp}[$V$ Controls Other Quantities]\label{assmp4: control of nabla V and norm}
        Assume there exists $w>0$ such that for every $x\in D$
    \begin{equation}\label{eq: nabla V bounded}
        \lvert \nabla V\rvert \leq c(1+V^w),\;\;\;
        \lvert x \rvert\leq c(V(x)+1).
    \end{equation}
\end{assmp}
\begin{assmp}[Monotonicity]\label{assmp5: monotonicity}
    There exists $\mu\in \mathbb{R}$ such that
    \begin{equation}
        \langle b(x)-b(y),x-y\rangle\leq \mu \lvert x-y\rvert^2.
    \end{equation}
\end{assmp}
About these assumptions, we make the following comments
\begin{itemize}
    \item We let $V$ take values in $[1,\infty)$ to avoid possible complications when it is raised to a negative power and arbitrarily close to $0$, as shall occur in the proofs.
    \item Assumption \ref{assmp1: Lipschitz like behaviour} is very similar to the `polynomial Lipschitz' condition used in the tamed schemes literature \cite{articletamed, articledif}, with a polynomial replaced by $V$. As we show in section \ref{subsec: superlinear}, we may indeed take $V$ to be a polynomial when our coefficients are superlinear and obey a coercivity condition. One notes as an immediate corollary that this implies $\lvert \nabla b\rvert \leq cV$.
    \item Similarly, Assumption \ref{assmp2: higher smoothness} is found in the tamed schemes literature, and is needed to obtain rate $1$ in convergence results. Similarly as above, this assumption implies $\lvert \nabla ^2 b\rvert\leq cV^q$.
    \item Assumption \ref{assmp3: moment assmp} is the technical core of the paper. The first condition, that $V$ explodes on $\partial D$, is key for proving \eqref{eq: main equation} has a strong solution. The subsequent part demonstrates that the true solution preserves `exponential $V$-moments' of a certain order, which allows us to deduce that $V(X_t)$ takes large values infrequently. Whilst it appears complicated, this condition holds for both superlinear and Lennard-Jones type coefficients. The bound on the rate of explosive trajectories lessens as $\alpha, \beta>0$ increase.
    \item Assumption \ref{assmp4: control of nabla V and norm} is a technical condition that we can quite easily bake into the design of $V$.
    \item Finally Assumption \ref{assmp5: monotonicity} is also common in the tamed schemes literature, as it allows us to bound the distance between the true solution and the scheme after we have corrected discretisation errors. We note that the strategy of \cite{articleboundedmeas, 10.1214/22-AAP1867, Dareiotis2026} does not rely on an assumption of this form, but relies on a more sophisticated PDE argument for convergence. It is possible these ideas could be extended to our setting.
\end{itemize}

As we show in section \ref{subsec: superlinear} and \ref{subsec: lennard jones} these assumptions cover both the case where $b$ is superlinear and coercive, as well as a class of coefficients which blow up to infinity on a compact set.

\subsection{Main Results}
As hinted previously, our main results shows that explosive trajectories are rare, by showing convergence conditional on a \textit{set of almost full measure}. In the following we denote by $E[X\mid A]$ the conditional expectation
\begin{equation}
    E[X\mid A]:=P(A)^{-1}\int_A X(\omega)P(d \omega).
\end{equation}
where $X$ is a random variable defined on a probability space $(\Omega, P, \mathcal{F})$ and $A\subset \Omega$ is a set of positive measure.

\begin{thm}\label{thm: main result}
Fix a time horizon $T>0$, and suppose all the assumptions of Section \ref{subsec: assumptions} hold. Then there exists a strong solution $(X_t)_{t\in [0,T]}$ to the dynamics \eqref{eq: main equation}. Furthermore for every $n\in \mathbb{N}$ there exists an event $A_n\subset \Omega$ satisfying 
\begin{equation}
    P(A_n)\geq 1-ce^{-ln^{\alpha r}},
\end{equation}
for $r=\frac{1}{3(\alpha+w+1)}$ and $l=\frac{1}{5}\wedge \beta$, such that the Euler scheme $(X^n_t)_{t\in [0,T]}$ defined in \eqref{eq: Euler scheme} satisfies
\begin{equation}
\sup_{t\in [0,T]}E[\lvert X^n_t-X_t\rvert^2 \mid A_n]\leq cn^{-2}.
\end{equation}

\end{thm}

We shall take $A_n$ to be the event $V(X^n_t)$ and $\lvert W_t-W_{\kappa_n(t)}\rvert$ are both almost surely bounded (relative to $n$) for $t\in [0,T]$. Specifically, $A_n:=\{\tau_n> T\}\subset\Omega$ for the stopping time $\tau_n$ given as 
\begin{equation}\label{eq: stopping time scheme}
    \tau_n:=\inf\{t>0\vert V(X^n_t)> n^{r}\;\text{or}\;\lvert W_t-W_{\kappa_n(t)}\rvert> n^{\alpha r/2-1/2}\}.
\end{equation}
This leads to the following corollaries which show converge to the entire law of $X_t$, in a strong and weak sense.
\begin{cor}[Strong Convergence]\label{cor: strong convergence}
     Let all the assumptions of Section \ref{subsec: assumptions} hold as above. Fix a time horizon $T>0$ and threshold parameters $C,h>0$. Consider the truncated approximation \eqref{eq: truncated approximation}. Then one has
     \begin{equation}
         \sup_{t\in [0,T]}E\lvert Y^n_t-X_t\rvert^2\leq cn^{-2}.
     \end{equation}
\end{cor}
We now prove weak convergence result in Wasserstein distance. We consider the law of $X^n_t$ conditioned on an event. For a set $A$ of positive measure, let us define the law of the conditional random variable $X\mid A$ as
\begin{equation}
    \mathcal{L}(X\mid A)(B):=\mathcal{L}(X)(B\cap A)/\mathcal{L}(X)(A),
\end{equation}
for all measurable sets $B$. Then one has

\begin{cor}[Weak Convergence]\label{cor: weak convergence}
    For the same assumptions and parameters as above, define the event
    \begin{equation}
        B_n:=\{\sup_{t\in [0,T]}\lvert X^n_{\kappa_n(t)}\rvert\leq Cn^h\}\subset \Omega.
    \end{equation}
    Then for the Wasserstein distance $\mathcal{W}_2$ as given in \eqref{eq: wasserstein} one has
    \begin{equation}
     \sup_{t\in [0,T]} \mathcal{W}_2(\mathcal{L}(X^n_t\mid B_n), \mathcal{L}(X_t))\leq cn^{-2}.
    \end{equation}
\end{cor}
Note that with access to samples $X(\omega_1),X(\omega_2),...$ of $X$, one may simulate from the law of $X\mid A$ by discarding draws for which $\omega_i\not \in A$.

\section{Analysis}
\subsection{Moment Bounds}
We begin with a technical lemma we shall use to prove exponential $V$-moments for both the scheme and the true solution. We consider here a semi-martingale with bounded drift and sufficiently bounded diffusion (up to a stopping time).
\begin{prop}\label{prop: control of moments}
Let $Y_t\geq 0$ be a positive adapted stochastic process satisfying for every $t\geq 0$
\begin{equation}
    Y_t=Y_0+\int^t_0 a_sds+\int^t_0 \langle b_s, dW_s\rangle ,\nonumber
\end{equation}
where $a_s$ is an adapted scalar valued process, the process $b_s$ is adapted and takes values in $\mathbb{R}^k$, one has $E[Y_0]<\infty$, and $W_s$ is a $\mathbb{R}^k$-valued Brownian motion. Suppose there exists a stopping time $\tau$, a constant $K>0$, and for every $R>0$ a constant $c_R>0$, such that the following holds almost surely
\begin{equation}\label{eq: condition on coef}
    a_{t\wedge \tau}\leq c,\;\;\; b_{t\wedge \tau}1_{\{Y_{t\wedge \tau}\leq R\}}\leq c_R.
\end{equation}
Then there exists a constant $\tilde{c}=\tilde{c}(K,T)>0$ such that for every $T\geq 0$ and $u>0$
\begin{equation}
    P(\sup_{t\in [0,T]} Y_{t\wedge \tau}\geq u)\leq \tilde{c}/u, \;\;\;
    E[\sup_{s\in [0,T]}\sqrt{Y_{t\wedge \tau}}]\leq \tilde{c},\nonumber
\end{equation}
\end{prop}
\begin{proof}
    In order to show that the stochastic integral vanishes, let us reason up to the stopping time
     \begin{equation}
        \tau_R:=\inf\{t\geq 0\vert Y_t\geq R\}.
    \end{equation}
    Then $a_{t\wedge \tau\wedge \tau_R}\leq c$ and $b_{t\wedge \tau\wedge \tau_R}\leq c_R$, so the stochastic integral up to $t\wedge \tau\wedge \tau_R$ is a true martingal, and therefore vanishes, and one has
    \begin{equation}
  \sup_{t\in [0,T]} E[Y_{t\wedge \tau \wedge\tau_R}]\leq K T,
\end{equation}
Using the continuity of sample paths and Markov's inequality, for $u\geq R$, it follows that
\begin{equation}
    P(\sup_{t\in [0,T]} Y_{t\wedge \tau}\geq u)= P( Y_{T\wedge \tau \wedge\tau_R}\geq u)\leq E[Y_{t\wedge \tau \wedge\tau_R}]/u=KT/u.
\end{equation}
Since we may take $R$ arbitrarily small and non-negative, this holds for all $u> 0$. For the second result note for $u\geq R$ as before
\begin{equation}
    P(\sup_{s\in [0,t]} \sqrt{Y_{s\wedge \tau}}\geq u)= P( Y_{t\wedge \tau \wedge\tau_R}\geq u^2)\leq KT/u^2.
\end{equation}
Applying this bound for $u\geq 1$, and bounding the LHS by $1$ for $u\in [0,1]$, one may integrate $u$ over $[0,\infty)$ to obtain the second bound by the continuous tail-sum formula.
\end{proof}
We may use this to show \eqref{eq: main equation} has a unique strong solution, and additionally has exponential $V$-moments.
\begin{lemma}[Moment Bounds for Continuous Process]\label{lemma: exponential moments for true solution}
Let Assumption \ref{assmp1: Lipschitz like behaviour}, \ref{assmp3: moment assmp} and \ref{assmp4: control of nabla V and norm} hold. Then \eqref{eq: main equation} has a unique strong solution for all $t\geq 0$ and furthermore
\begin{equation}
    P(\sup_{t\in [0,T]} e^{\beta V(X_t)^\alpha}\geq u)\leq c/u, \;\;\;
    E[\sup_{s\in [0,T]}e^{\frac{\beta }{2}V(X_t)^\alpha}]\leq c,\nonumber
\end{equation}
\end{lemma}
\begin{proof}
Let $x,y\in D$ be arbitrary. Then, by Assumption \ref{assmp1: Lipschitz like behaviour} one has
\begin{equation}
    \lvert b(x)\rvert\leq \lvert b(y)\rvert+\lvert b(x)-b(y)\rvert\leq \lvert b(y)\rvert+c(V(x)+V(y))\lvert x-y\rvert.
\end{equation}
Therefore, fixing $y\in D$, by the triangle inequality one has the following bound that is uniform over $x\in D$
\begin{equation}\label{eq: key b bound}
    \lvert b(x)\rvert\leq c(1+V(x))(1+\lvert x \rvert).
\end{equation}
Therefore, for every $R>0$ one has
\begin{equation}
    \lvert b(x)\rvert 1_{\{V(x)\leq R\}}\leq c(1+R)(1+\lvert x \rvert).
\end{equation}
Let us assume that the initial condition $X_0\in D$ is a constant. Then applying the result of \cite{ZHANG20051805}, one sees that for every $R>0$ the dynamics
\begin{equation}\label{eq: truncated dynamics}
        d\tilde{X}^R_t=b(\tilde{X}^R_t)1_{\{\lvert b(X_t)\rvert\leq R\}}dt+dW_t,\;\;\;\tilde{X}^R_0=X_0,
    \end{equation}
    have a unique strong solution on $\mathbb{R}^d$. Therefore, we may use the solution to the dynamics \eqref{eq: truncated dynamics} to construct a unique solution to \eqref{eq: main equation} up to the stopping time $\tilde{ \tau}_R$, where 
    \begin{equation}
        \tilde{\tau}_R=\inf\{t\geq 0\vert \;\lvert V(X_t)\rvert\geq R\}.
    \end{equation}
    Note we do not require $X_0$ to be constant to perform this, since the strong solution holds for arbitrary $V(X_0)\leq R$, and if $V(X_0(\omega))\geq R$ then $\tilde{\tau}_R(\omega)=0$. 
Therefore, the following is well-defined
    \begin{align}
        \exp(\beta V(X_{t\wedge \tilde{\tau}_R})^\alpha)&=\exp(\beta V(X_0)^\alpha)+\int^{t\wedge \tilde{\tau}_R}_0 (\mathcal{L}\exp(\beta V^\alpha))(X_t)ds\nonumber \\
        &+\int^{t\wedge \tilde{\tau}_R}_0 \beta \alpha V(X_s)^{\alpha-1}\exp(\beta V(X_s)^\alpha) \langle \nabla V(X_s), dW_s\rangle.
    \end{align}
    Note that the first integrand is bounded independently of $R>0$ by Assumption \ref{assmp3: moment assmp}, and the second is bounded as a function of $R$ by the definition of the stopping time and Assumption \ref{assmp4: control of nabla V and norm}. Therefore, we may Proposition \ref{prop: control of moments} to obtain
    \begin{equation}\label{eq: bound on cont stopped}
    P(\sup_{t\in [0,T]} \exp(\beta V(X_{t\wedge \tilde{\tau}_R})^\alpha)\geq u)\leq c/u,
    \end{equation}
    \begin{equation}
    E[\sup_{s\in [0,T]}\exp(\frac{\beta}{2} V(X_{t\wedge \tilde{\tau}_R})^\alpha)]\leq c,
\end{equation}
where the constant here does not depend on $R>0$. Therefore
\begin{align}
    P(\tilde{\tau}_R\in [0,T])&=P(\sup_{t\in [0,T]}V(X_{t\wedge \tilde{\tau}_R})\geq R)\nonumber \\
    &=P(\sup_{t\in [0,T]} \exp(\beta V(X_{t\wedge \tilde{\tau}_R})^\alpha\geq e^{\beta R^\alpha})\nonumber \\
    &\leq ce^{-\beta R^\alpha}.
\end{align}
Then
\begin{equation}
P(\sup_{R>0}\tilde{\tau}_R\in [0,T])=P(\cap_{R>0}\{V(X_{t\wedge \tilde{\tau}_R})\geq R\})=0,
\end{equation}
so a unique strong solution to \eqref{eq: main equation} can be given by setting $X_t=\lim_{R\to\infty}X_{t\wedge \tilde{\tau}_R}$. The stated bounds then follow from Fatou's lemma, writing the probability in \eqref{eq: bound on cont stopped} as the expectation of an indicator function.
\end{proof}

\begin{lemma}[Moment Bounds for the Stopped Scheme]\label{lem: moment bounds scheme}
    Let assumptions \ref{assmp1: Lipschitz like behaviour}, \ref{assmp2: higher smoothness}, \ref{assmp3: moment assmp} and \ref{assmp4: control of nabla V and norm} hold. Then
    \begin{equation}
    P(\sup_{t\in [0,T]} e^{\beta V(X^n_{t\wedge \tau_n})^\alpha}\geq u)\leq c/u, \;\;\;
    E[\sup_{s\in [0,T]}e^{\frac{\beta }{2}V(X^n_{t\wedge \tau_n})^\alpha}]\leq c.
\end{equation}
\end{lemma}
\begin{proof}
Using the stopping time $\tau_n$ given above in \eqref{eq: stopping time scheme}, one may correct the discretisation lag to write
 \begin{align}\label{eq: discrete exp}
        &e^{\beta V(X^n_{t\wedge \tau_n})^\alpha}=e^{\beta V(X_0)^\alpha}+\int^{t\wedge \tau_n}_0(\mathcal{L}e^{\beta V^\alpha})(X^n_s)ds\nonumber \\
        &+\int^{t\wedge \tau_n}_0 \alpha V(X^n_s)^{\alpha-1} e^{\beta V(X^n_s)^\alpha}\langle \nabla V(X^n_s), b(X^n_{\kappa_n(s)})-b(X^n_s)\rangle ds \nonumber \\
         &+\int^{t\wedge \tau_n}_0\alpha  V(X^n_s)^{\alpha-1} e^{\beta V(X^n_s)^\alpha}\langle \nabla V(X^n_s), dW_s\rangle.
    \end{align}
    Using Assumption \ref{assmp1: Lipschitz like behaviour}, and expanding the increment, one has
    \begin{equation}
        \lvert b(X^n_s)-b(X^n_{\kappa_n(s)}\rvert\leq c(V(X^n_s)+V(X^n_{\kappa_n(s)}))(n^{-1}\lvert b(X^n_{\kappa_n(t)})\rvert+\lvert W_s-W_{\kappa_n(t)}\rvert).\nonumber
    \end{equation}
    Using the bounds on $\nabla V$ and $\lvert x\rvert$ in Assumption \ref{assmp4: control of nabla V and norm}, bounding $V^{\alpha-1}\leq V^\alpha$ by the assumption $V\geq 1$, applying \eqref{eq: key b bound}, and finally applying Young's inequality, one has
    \begin{align}
&V(X_s)^{\alpha-1} \langle \nabla V(X^n_s), b(X^n_{\kappa_n(s)})-b(X^n_s)\rangle \leq V(X_s)^{\alpha+w}\lvert b(X^n_{\kappa_n(s)})-b(X^n_s)\rvert\nonumber \\
&\leq n^{-1}(V(X_s)^{\alpha+w+3}+V(X_{\kappa_n(s)})^{\alpha+w+3})\nonumber \\
&+\lvert W_s-W_{\kappa_n(s)}
\rvert (V(X_s)^{\alpha+w+1}+V(X_{\kappa_n(s)})^{\alpha+w+1}).
    \end{align}
    Then using the definition of the stopping time \eqref{eq: stopping time scheme} one sees that the expression above is uniformly bounded up to $\tau_n$, and therefore
     \begin{align}\label{eq: discrete exp 2}
        &e^{\beta V(X^n_{t\wedge \tau_n})^\alpha}=e^{\beta V(X_0)^\alpha}+\int^{t\wedge \tau_n}_0((\mathcal{L}e^{\beta V^\alpha})(X^n_s)+c e^{\beta V(X^n_s)^\alpha})ds\nonumber \\
         &+\int^{t\wedge \tau_n}_0\alpha  V(X^n_s)^{\alpha-1} e^{\beta V(X^n_s)^\alpha}\langle \nabla V(X^n_s), dW_s\rangle,
    \end{align}
so that by Assumption \ref{assmp3: moment assmp} one has
\begin{align}\label{eq: discrete exp 3}
        &e^{\beta V(X^n_{t\wedge \tau_n})^\alpha}=e^{\beta V(X_0)^\alpha}+c\int^{t\wedge \tau_n}_0ds\nonumber \\
         &+\int^{t\wedge \tau_n}_0\alpha  V(X^n_s)^{\alpha-1} e^{\beta V(X^n_s)^\alpha}\langle \nabla V(X^n_s), dW_s\rangle.
    \end{align}
    Therefore the assumptions of Proposition \ref{prop: control of moments} are fulfilled (using Assumption \ref{assmp4: control of nabla V and norm} for the drift coefficient), so the result follows.
\end{proof}
\begin{lemma}\label{lem: poly V moments}
For every $p>0$ and $T\in [0,T]$ one has
    \begin{equation}
       E[\sup_{t\in [0,T]}V(X_t)^p]\leq c,\;\;\; E[\sup_{t\in [0,T]}V(X^n_{t\wedge \tau_n})^p]\leq c,
    \end{equation}

    \begin{equation}
   E[\sup_{t\in [0,T]}\lvert X_t\rvert^p]\leq c,\;\;\; E[\sup_{t\in [0,T]}\lvert X^n_{t\wedge \tau_n}\rvert^p]\leq c.
    \end{equation}
\end{lemma}
\begin{proof}
    The result follows using the results of this section, Assumption \ref{assmp4: control of nabla V and norm} and the fact $V(x) ^p\leq c_pe^{\beta V(x)^\alpha}$.
\end{proof}
\subsection{Stopping Analysis}\label{subsec: stopping analysis}
In this section we bound the probability the stopping time $\tau_n$, introduced in \eqref{eq: stopping time scheme}, lies in the interval $[0,T]$.
\begin{prop}\label{prop: stopping time unlikely}
    Let assumptions \ref{assmp1: Lipschitz like behaviour}, \ref{assmp2: higher smoothness}, \ref{assmp3: moment assmp} and \ref{assmp4: control of nabla V and norm} hold. Then for the stopping time \eqref{eq: stopping time scheme} one has
    \begin{equation}
        P(\tau_n\in [0,T])\leq c\exp (-ln^{\alpha r}),
    \end{equation}
    for $r,l>0$ as in Theorem \ref{thm: main result}.
\end{prop}
\begin{proof}
One notes, if the stopping time $\tau_n$ satisfies $\tau_n\in [0,T]$, the conditions for stopping must be fulfilled before stopping, that is on $[0,\tau_n]$. Therefore
\begin{align}
 &\{\tau_n\in [0,T]\}= \{\sup_{t\in [0,T\wedge \tau_n]}V(X^n_t)\}\cup\{\sup_{t\in [0,T\wedge \tau_n]}\lvert W_t-W_{\kappa_n(t)}\rvert\geq n^{\alpha r/2-1/2}\}\nonumber \\
 &= \{\sup_{t\in [0,T]}V(X^n_{t\wedge \tau_n})\}\cup\{\sup_{t\in [0,T]}\lvert W_{t\wedge \tau_n}-W_{\kappa_n(t)\wedge \tau_n}\rvert\geq n^{\alpha r/2-1/2}\},
\end{align}
and so by the union bound one has
\begin{align}\label{eq: splitting of bound}
  P(\tau_n\leq T)&\leq P(\sup_{t\in [0,T]}V(X^n_{t\wedge \tau_n})\geq n^r)\nonumber \\
  &+P(      \sup_{t\in [0,T]}\lvert W_{t\wedge \tau_n}-W_{\kappa_n(t)\wedge \tau_n}\rvert\geq n^{\alpha r/2-1/2}).
    \end{align}
    For the first term of \eqref{eq: splitting of bound}, by Lemma \ref{lem: moment bounds scheme} and Markov's inequality one has
    \begin{align}\label{eq: first bound stopping}
    P(\sup_{t\in [0,T]}V(X^n_{t\wedge \tau_n})\geq n^r)&=P(\sup_{t\in [0,T]}e^{\beta V(X^n_{t\wedge \tau_n})^\alpha}\geq e^{\beta n^{\alpha r}})\nonumber \\
    &\leq ce^{-\beta n^{\alpha r}}.
    \end{align}
     For the second term of \eqref{eq: splitting of bound}, noting that $[0,T]$ contains no more than $T(n+1)$ intervals of length $n^{-1}$ one has
    \begin{align}\label{eq: bound on stopping second term}
        P(\sup_{t\in [0,T]}\lvert W_{t\wedge \tau_n}-&W_{\kappa_n(t)\wedge \tau_n}\rvert\geq n^{r/2-1/2})\leq P(\sup_{t\in [0,T]}\lvert W_t-W_{\kappa_n(t)}\rvert\geq n^{r/2-1/2}) \nonumber \\
       &\leq cnP(\sup_{t\in [0,n^{-1}]}\lvert W_t\rvert\geq n^{r/2-1/2})\nonumber \\
        &\leq  cnP(\sup_{t\in [0,n^{-1}]}\lvert W_t\rvert^2\geq n^{r-1})\nonumber \\
        &\leq cn P(e^{n \sup_{t\in [0,n^{-1}]}\lvert W_t\rvert^2/4})\geq e^{n^r/4}).
    \end{align}
    Now let us show that 
    \begin{equation}\label{eq: exp integrability}
    E[\exp(n^{1/2}\sup_{t\in [0,n^{-1}]}\lvert W_t\rvert^2/4)]\leq c,
    \end{equation}
    so as to apply Markov's inequality in the expression above. One has
    \begin{equation}\label{eq: product identitity}
\exp(n\sup_{t\in [0,n^{-1}]}\lvert W_t\rvert^2/4)= \prod_{i=1}^d\exp(n\sup_{t\in [0,n^{-1}]}\lvert W^i_t\rvert^2/4),
    \end{equation}
    so by the independence of each coordinate process of $W_t$, and since we do not track dimensional dependence in the generic constant $c>0$, it is sufficient to show
    \begin{equation}
       E[\exp(n\sup_{t\in [0,n^{-1}]}\lvert W^i_t\rvert^2/4)]\leq c.
    \end{equation}
    By the continuous tail-sum formula one has
    \begin{align}
       &E[\exp(n\sup_{t\in [0,n^{-1}]}\lvert W^i_t\rvert^2/4)]=\int^\infty_0P(\exp(n\sup_{t\in [0,n^{-1}]}\lvert W^i_t\rvert^2/4)> a)da\nonumber \\
       &=\int^\infty_0P(\sup_{t\in [0,n^{-1}]}\lvert W^i_t\rvert> 2\sqrt{\log(a)n^{-1}})da\nonumber \\
       &\leq\int^\infty_0(P(\sup_{t\in [0,n^{-1}]} W^i_t> 2\sqrt{\log(a)n^{-1}})+P(\sup_{t\in [0,n^{-1}]}(- W^i_t)> 2\sqrt{\log(a)n^{-1}})da\nonumber
    \end{align}
Therefore by the reflection principle, and since Brownian motion is equal in law to its negative writing $\xi\sim \mathcal{N}(0,1)$ for a standard Gaussian one has
\begin{align}
       E[\exp(n\sup_{t\in [0,n^{-1}]}\lvert W^i_t\rvert^2/4)]&\leq 4\int^\infty_0P(W^i_{n^{1/2}}> 2\sqrt{\log(a)n^{-1}})da\nonumber\\
       &\leq 4\int^\infty_0P(\xi> 2\sqrt{\log(a)})da.
    \end{align}
Furthermore, one may bound the integrand by $1$ for $a\in [0,1]$, and use the fact $E[e^{\xi^2/3}]\leq c$ to apply Markov's inequality on $a\in (1,\infty)$, to obtain
\begin{align}
       E[\exp(n\sup_{t\in [0,n^{-1}]}\lvert W^i_t\rvert^2/4)]&\leq 4+4\int^\infty_1P(e^{\xi^2/3}> a^{4/3})da\nonumber \\
       &\leq 4+c\int_1^\infty a^{-4/3}da\nonumber \\
       &\leq c.
    \end{align}
Therefore by \eqref{eq: product identitity} we have shown \eqref{eq: exp integrability}, and so by Markov's inequality one uses \eqref{eq: bound on stopping second term} to bound
\begin{equation}\label{eq: second bound stopping}
    P(\sup_{t\in [0,T]}\lvert W_{t\wedge \tau_n}-W_{\kappa_n(t)\wedge \tau_n}\rvert\geq n^{r/2-1/2})\leq cne^{-n^r/4}\leq ce^{-n^r/5}.
\end{equation}
Combining \eqref{eq: first bound stopping} and \eqref{eq: second bound stopping} into \eqref{eq: splitting of bound} the result follows.
\end{proof}
\subsection{Convergence Analysis}\label{subsec: convergence analysis}
Here we show convergence up to the stopping time $\tau_n$, which we have shown in Proposition \ref{prop: stopping time unlikely} is unlikely to lie in $[0,T]$ for any fixed $T>0$. The following argument is unfortunately quite technical, although the ideas are standard in the numerics literature.
\begin{prop}\label{prop: convergence up to stopping}
    Let assumptions \ref{assmp1: Lipschitz like behaviour}, \ref{assmp2: higher smoothness}, \ref{assmp3: moment assmp}, \ref{assmp4: control of nabla V and norm} and \ref{assmp5: monotonicity} hold. Then for every $t\in [0,T]$ and $n\in \mathbb{N}$ one has
    \begin{equation}
        \sup_{t\in [0,T]}E\lvert X^n_{t\wedge \tau_n}-X_{t\wedge \tau_n}\rvert^2\leq cn^{-2}.
    \end{equation}
\end{prop}
\begin{proof}
    To simplify notation let us denote
    \begin{equation}
        e_t:=X^n_t-X_t.
    \end{equation}
    By Itô's formula one has
    \begin{align}
   \lvert e_{t\wedge \tau_n}\rvert^2&=\int^{t\wedge\tau_n}_02\langle b(X^n_s)-b(X_s), e_s\rangle ds\nonumber \\
    &+\int^{t\wedge \tau_n}_02\langle b(X^n_{\kappa_n(s)})-b(X^n_s), e_s\rangle ds.
\end{align}
Therefore, applying Assumption \ref{assmp5: monotonicity} one has
 \begin{align}
   \lvert e_{t\wedge \tau_n}\rvert^2&=\int^{t\wedge\tau_n}_0c\lvert e_s\rvert^2  ds+\int^{t\wedge \tau_n}_02\langle b(X^n_{\kappa_n(s)})-b(X^n_s), e_s\rangle ds.
\end{align}
The rest of the proof shall involve exploiting subtle probabilistic cancellations in the second term. To this end, let us begin by splitting as
\begin{align}\label{eq: splitting 1}
   & \langle e_s, b(X_s)-b(X_{\kappa_n(s)})\rangle\nonumber \\
   &\leq \langle e_s, b(X_s)-b(X_{\kappa_n(s)})-\nabla b(X_{\kappa_n(s)})(X_s-X_{\kappa_n(s)})\rangle\nonumber \\
   &+\langle e_s, \nabla b(X_{\kappa_n(s)})(X_s-X_{\kappa_n(s)})\rangle.
\end{align}
We now split once more by writing the second term as
\begin{align}\label{eq: splitting 2}
    &\langle e_s, \nabla b(X_{\kappa_n(s)})(X_s-X_{\kappa_n(s)})\rangle\nonumber \\
    &=\biggr \langle e_{\kappa_n(s)}+\int^s_{\kappa_n(s)}(b(X^n_{\kappa_n(u)})-b(X_{\kappa_n(u)}) )du+\int^s_{\kappa_n(s)}(b(X_{\kappa_n(u)})-b(X_u) )du,\nonumber \\
    &\nabla b(X_{\kappa_n(s)})\int^s_{\kappa_n(s)}b(X_u)du+\nabla b(X_{\kappa_n(s)})(W_s-W_{\kappa_n(s)})\biggr \rangle.
\end{align}
This splitting of the inner product into two terms on the left and three on the right results in six terms. Along with the first term of the splitting \eqref{eq: splitting 1} this yields seven terms. Of these seven terms, let us firstly bound two of them by Cauchy-Schwarz's inequality and Young's inequality. Specifically, one bounds the first term in the splitting \eqref{eq: splitting 1} as
\begin{align}
    &\langle e_s, b(X_s)-b(X_{\kappa_n(s)})-\nabla b(X_{\kappa_n(s)})(X_s-X_{\kappa_n(s)})\rangle\nonumber \\
    &\leq c\lvert e_s\rvert^2+c\lvert b(X_s)-b(X_{\kappa_n(s)})-\nabla b(X_{\kappa_n(s)})(X_s-X_{\kappa_n(s)}) \rvert^2,
\end{align}
and the product of the first term on the left side of the product, and first on the right side in \eqref{eq: splitting 2}, to obtain
\begin{align}
    &\langle e_{\kappa_n(s)}, \nabla b(X_{\kappa_n(s)})\int^s_{\kappa_n(s)}b(X_u)du\rangle\nonumber \\
    &\leq c\lvert e_{\kappa_n(s)}\rvert^2+c\biggr \lvert \nabla b(X_{\kappa_n(s)})\int^s_{\kappa_n(s)}b(X_u)du\biggr \rvert^2.
\end{align}
As a result of these calculations, one may write
\begin{align}
   & \langle e_s, b(X_s)-b(X_{\kappa_n(s)})\rangle\nonumber \\
   &\leq c\lvert e_s\rvert^2+c\lvert e_{\kappa_n(s)}\rvert^2\nonumber \\
   &+I^1_s+I^2_s,
\end{align}
for
\begin{align}
    I^1_s&:=\langle e_{\kappa_n(s)}, \nabla b(X_{\kappa_n(s)})(W_s-W_{\kappa_n(s)})\rangle\nonumber \\
    &+\langle \int^s_{\kappa_n(s)}(b(X^n_{\kappa_n(u)})-b(X_{\kappa_n(u)}) )du, \nabla b(X_{\kappa_n(s)})(W_s-W_{\kappa_n(s)})\rangle\nonumber \\
\end{align}
and
\begin{align}\label{eq: I^2_s defn}
    I^2_s&:=c\lvert b(X_s)-b(X_{\kappa_n(s)})-\nabla b(X_{\kappa_n(s)})(X_s-X_{\kappa_n(s)}) \rvert^2\nonumber \\
    &+c\biggr \lvert \nabla b(X_{\kappa_n(s)})\int^s_{\kappa_n(s)}b(X_u)du\biggr \rvert^2\nonumber \\
    &+\langle \int^s_{\kappa_n(s)}(b(X_{\kappa_n(u)})-b(X_u) )du,\nabla b(X_{\kappa_n(s)})(W_s-W_{\kappa_n(s)})\rangle\nonumber \\
    &+\langle \int^s_{\kappa_n(s)}(b(X^n_{\kappa_n(u)})-b(X_{\kappa_n(u)}) )du, \nabla b(X_{\kappa_n(s)})\int^s_{\kappa_n(s)}b(X_u)du\rangle\nonumber \\
    &\langle \int^s_{\kappa_n(s)}(b(X_{\kappa_n(u)})-b(X_u) )du, \nabla b(X_{\kappa_n(s)})\int^s_{\kappa_n(s)}b(X_u)du\rangle\rangle .
\end{align}
Note there are seven error terms as required, as well as $c\lvert e_s\rvert^2+c\lvert e_{\kappa_n(s)}\rvert^2$. Using these expressions one may write
\begin{align}\label{eq: bound on stopped difference}
   \lvert e_{t\wedge \tau_n}\rvert^2&=\int^{t\wedge\tau_n}_0c(\lvert e_s\rvert^2 +\lvert e_{\kappa_n(s)}\rvert^2)ds\nonumber \\
    &+\int^{t}_0 I^1_{s\wedge \tau_n}ds-\int^t_{t\wedge\tau_n}I^1_{s\wedge\tau_n}ds+\int^{t\wedge \tau_n}_0I^2_{s\wedge \tau_n}ds.
\end{align}
Since $E[W_{s\wedge \tau_n}-W_{\kappa_n(s)\wedge \tau_n}\mid\mathcal{F}_{\kappa_n(s)}]=0$, one has that $\int^{t}_0 I^1_{s\wedge \tau_n}ds=0$. Furthermore, note that by Assumption \ref{assmp1: Lipschitz like behaviour} one has
\begin{equation}\label{eq: nabla b identity}
    \lvert \nabla b(x)\rvert \leq V(x),
\end{equation}
so that using this as well as Lemma \ref{lem: poly V moments} and standard inequalities one has $\sup_{s\in [0,T]}E[\lvert I^1_{t\wedge \tau_n}\rvert^2]\leq c$. Therefore by Hölder's inequality and Proposition \ref{prop: stopping time unlikely}, and since exponential decay dominates polynomial decay, one has
\begin{align}
   E\biggr [\int^t_{t\wedge\tau_n}I^1_{s\wedge\tau_n}\biggr]&= E\biggr [1_{\tau_n\leq T}\cdot\int^t_{t\wedge\tau_n}I^1_{s\wedge\tau_n}\biggr]\nonumber \\
   &\leq c\sqrt{P(\tau_n\leq T)\int^t_0 E[\lvert I^1_{s\wedge \tau_n}\rvert^2]ds}\nonumber \\
   &\leq cn^{-2}.
\end{align}
Finally let us show $\sup_{s\in [0,T]}E[\lvert I^2_{t\wedge \tau_n}\rvert^2]\leq cn^{-2}$. This follows from Lemma \ref{lem: poly V moments}, the assumptions, the identity \eqref{eq: nabla b identity} and standard inequalities. For instance, one bounds the third term of $I^2_{s\wedge \tau_n}$ to the appropriate order in the following way
\begin{align}
    &\lvert \langle \int^{s\wedge \tau_n}_{\kappa_n(s\wedge \tau_n)}(b(X_{\kappa_n(u)})-b(X_u) )du,\nabla b(X_{\kappa_n(s)})(W_s-W_{\kappa_n(s)})\rangle\rvert \nonumber \\
    &\leq \int^{s}_{\kappa_n(s)}(V(X_{\kappa_n(u)\wedge \tau_n})+V(X_{u\wedge \tau_n}))  \times \nonumber \\
   & +\biggr ( n^{-1}\lvert b(X_{\kappa_n(u)\wedge \tau_n})\rvert+\lvert W_{s\wedge \tau_n}-W_{\kappa_n(s)\wedge \tau_n}\rvert \biggr )\times \nonumber \\
    &V(X_{\kappa_n(s)\wedge \tau_n})\lvert W_{s\wedge \tau_n}-W_{\kappa_n(s)\wedge \tau_n}\rvert du.
\end{align}
One sees that the expression inside the integral is bounded in expectation by $cn^{-1}$, and since the integral is over an interval of length $n^{-1}$ the appropriate bound holds. Similar results hold for each term of $I^2_s$. Putting this together one has from \eqref{eq: bound on stopped difference}, uniformly in $t\in [0,T]$, the following bound
\begin{align}
   \lvert e_{t\wedge \tau_n}\rvert^2&=\int^{t\wedge\tau_n}_0c(\lvert e_s\rvert^2 +\lvert e_{\kappa_n(s)}\rvert^2)ds+cn^{-2}\nonumber \\
   &\leq c\int^t_0\sup_{u\in [0,s]}\lvert e_{u\wedge \tau_n}\rvert^2 ds+cn^{-2}.
\end{align}
Therefore by Grönwall's inequality the result follows.
\end{proof}
\subsection{Proofs of Main Results}\label{subsec: proofs}
\subsubsection{Proof of Theorem \ref{thm: main result}}
First note by Proposition \ref{prop: stopping time unlikely} one has that the set $A_n=\{\tau_n> T\}\subset \Omega$ given in Theorem \ref{thm: main result} is of positive measure. Therefore
\begin{align}
 \sup_{t\in [0,T]}E[\lvert X^n_t-X_t\rvert^2 \vert A_n]&=P(A_n)^{-1}\sup_{t\in [0,T]}\int_{A_n}\lvert X^n_t(\omega)-X_t(\omega)\rvert^2dP(\omega)\nonumber \\
 &= P(A_n)^{-1}\sup_{t\in [0,T]}E[\lvert X^n_t-X_t\rvert^21_{A_n}]\nonumber \\
&= P(A_n)^{-1}\sup_{t\in [0,T]}E[\lvert X^n_{t\wedge \tau_n}-X_{t\wedge \tau_n}\rvert^21_{A_n}]\nonumber \\
 &\leq P(A_n)^{-1}\sup_{t\in [0,T]}E[\lvert X^n_{t\wedge \tau_n}-X_{t\wedge \tau_n}\rvert^2].
\end{align}
Therefore, combining Proposition \ref{prop: stopping time unlikely} and Proposition \ref{prop: convergence up to stopping} the result follows for all sufficiently large $n$, since $P(A_n)\to 1$ as $n\to \infty$. Furthermore, since $\sup_{n\geq 1}E[\sup_{t\in [0,T]}\lvert X^n_t\rvert^2]<\infty$ and $E[\sup_{t\in [0,T]}\lvert X_t\rvert^2]<\infty$ by Lemma \ref{lem: poly V moments}, the result can be extended to all $n\in \mathbb{N}$ by choosing the constant $c>0$ in the Theorem sufficiently large.

\subsubsection{Proof of Corollary \ref{cor: strong convergence}}
For $C,r>0$ as given in the corollary statement, as in Corollary \ref{cor: weak convergence} let us set
\begin{equation}
B_n:=\{\sup_{t\in [0,T]}\lvert X^n_{\kappa_n(t)}\rvert \leq Cn^r\}\subset \Omega,
\end{equation}
so that $Y^n_t=X^n_t\cdot 1_{B_n}$ and one has
\begin{align}
\sup_{t\in [0,T]}E[\lvert Y^n_t-X_t\rvert^2 ]&\leq c\sup_{t\in [0,T]}E[\lvert X^n_t\cdot 1_{B_n}-X^n_t\rvert^2 1_{A_n\cap B_n}] \nonumber \\
&+ c \sup_{t\in [0,T]}E[\lvert X^n_t\cdot 1_{B_n}-X_t\rvert^21_{(A_n\cap B_n)^c}].
\end{align}
For the first term, using the definition of $A_n$ and Proposition \ref{prop: convergence up to stopping} one has
\begin{align}\label{eq: splitting of difference}
    \sup_{t\in [0,T]}E[\lvert X^n_t\cdot 1_{B_n}-X^n_t\rvert^2 1_{A_n\cap B_n}] &=\sup_{t\in [0,T]}E[\lvert X^n_{t\wedge \tau_n}-X^n_{t\wedge \tau_n}\rvert^2 1_{A_n\cap B_n}]\nonumber \\
    &\leq \sup_{t\in [0,T]}E[\lvert X^n_{t\wedge \tau_n}-X^n_{t\wedge \tau_n}\rvert^2 ]\nonumber \\
    &\leq cn^{-2}.
\end{align}
For the second term one applies Hölder's inequality to obtain
\begin{align}\label{eq: splitting}
    &\sup_{t\in [0,T]}E[\lvert X^n_t\cdot 1_{B_n}-X_t\rvert^21_{(A_n\cap B_n)^c}]\nonumber \\
    &\leq  \sup_{t\in [0,T]}\sqrt{E[\lvert X^n_t\cdot 1_{B_n}-X_t\rvert^4]P((A_n\cap B_n)^c)}.\nonumber
\end{align}
Therefore, using the bound on $X^n_t\cdot 1_{B_n}$ implied by construction and Lemma \ref{lem: poly V moments} 
\begin{equation}\label{eq: first factor}
\sup_{t\in [0,T]}\sqrt{E[\lvert X^n_t\cdot 1_{B_n}-X_t\rvert^4]}\leq c(1+C^2n^{2r}).
\end{equation}
Furthermore, one has
\begin{equation}
    P((A_n\cap B_n)^c)\leq P(A_n^c)+P(B_n^c).
\end{equation}
and by Markov's inequality and Lemma \ref{lem: poly V moments}
\begin{align}
    P(B_n^c)&=P(\sup_{t\in [0,T]}\lvert X^n_{\kappa_n(t)}\rvert > Cn^r)\nonumber \\
    &=P(\sup_{t\in [0,T]}\lvert X^n_{\kappa_n(t)}\rvert^{2/r} > C^{2/r}n^{2})\nonumber \\
    &\leq cn^{-2}.
\end{align}
Therefore by the bound on $A_n$ from Theorem \ref{thm: main result} one has for every $a>0$
\begin{equation}\label{eq: second factor}
    P((A_n\cap B_n)^c)\leq cn^{-2}.
\end{equation}
Combining \eqref{eq: first factor} and \eqref{eq: second factor} one has
\begin{equation}\label{eq: bound on second term}
\sup_{t\in [0,T]}E[\lvert X^n_t\cdot 1_{B_n}-X_t\rvert^21_{(A_n\cap B_n)^c}]\leq cn^{-2},
\end{equation}
at which point the result follows in combination with \eqref{eq: splitting of difference}.

\subsection{Proof of Corollary \ref{cor: weak convergence}}
To show convergence in Wasserstein distance, for every $t\in [0,T]$ we must find a coupling between $\mathcal{L}(X^n_t\vert B_n)$ and $\mathcal{L}(X_t)$. Let us build this coupling on the probability space $(\Omega, P, \mathcal{F})$ that we have constructed the continuous time processes. Let $Z_t\sim \mathcal{L}(X^n_t\vert B_n)$ be defined as
\begin{itemize}
    \item $Z_t(\omega)=X^n_t(\omega)$ when $\omega \in B_n$
    \item For every $\omega \not \in B_n$ set $Z_t=Q_t$, where $Q_t$ is a realisation of $\mathcal{L}(X^n_t\vert B_n)$ independent of $X_t$ and $X^n_t$
\end{itemize}
One notes that $Z_t$ does indeed have law $\mathcal{L}(X^n_t\vert B_n)$, since $Z_t \vert B_n \sim \mathcal{L}(X^n_t\vert B_n)$ and $Z_t \vert B_n^c \sim \mathcal{L}(X^n_t\vert B_n)$. Furthermore
\begin{align}
\sup_{t\in [0,T]}W_2(\mathcal{L}(X^n_t\vert B_n),X_t)^2   &\leq \sup_{t\in [0,T]}E[\lvert Z_t-X_t\rvert^2]\nonumber \\
&\leq \sup_{t\in [0,T]}E[\lvert X^n_t-X_t\rvert^21_{A_n\cap B_n}]\nonumber \\
&+\sup_{t\in [0,T]}E[\lvert Z_t-X_t\rvert^21_{(A_n\cap B_n)^c}].
\end{align}
The result then follows from an almost identical argument to the one presented previously: the first term follows from \eqref{eq: splitting of difference}, and the second follows from almost the same argument as used to obtain \eqref{eq: bound on second term}, since one has $\lvert Z_t\rvert \leq Cn^r$ almost surely by construction.

\section{Examples}\label{sec: examples}
We consider two classes of well-studied non-Lipschitz drift coefficients, whose non-Lipschitz behaviour leads to the explosion of Euler scheme discretisations.
\subsection{Superlinear Coefficients}\label{subsec: superlinear}
We begin by considering coefficients that grow superlinearly, for example
\begin{equation}\label{eq: example SDE superlinear}
    dX_t=-X_t^3+dW_t.
\end{equation}
The instability of Euler approximations of dynamics of this kind is shown in \cite{ae1c5b69-79b7-3bf4-87cb-ba6dd72c5ca2}. Specifically, we assume
\begin{Aassumption}[Polynomial Lipschitz]\label{A-assmp: polynomial lipschitz}
    The domain $D=\mathbb{R}^d$, and there exists $\eta>0$ such that and every $x,y \in \mathbb{R}^d=D$
\begin{equation}
    \lvert b(x)-b(y)\rvert \leq c(1+\lvert x \rvert^\eta+\lvert y \rvert^\eta)\lvert x-y\rvert.
\end{equation}
\end{Aassumption}
\begin{Aassumption}[Polynomial Smoothness]\label{A-assmp: polynomial lipschitz smoothness}
    There exists $\xi\geq \max\{\eta, 2\}$ such that the following smoothness condition holds for every $x,y \in D$
\begin{equation}
    \lvert \nabla b(x)-\nabla b(y)\rvert \leq c(1+\lvert x \rvert^\xi+\lvert y \rvert^\xi)\lvert x-y\rvert.
\end{equation}
\end{Aassumption}
\begin{Aassumption}[Coercivitity]\label{A-assmp: coercivity}
    There exists $L>0$ such that for all $x\in \mathbb{R}^d$
\begin{equation}
    \langle b(x),x\rangle \leq c-L\lvert x \rvert^2.
\end{equation}
\end{Aassumption}
The A-Assumptions \ref{A-assmp: polynomial lipschitz} and \ref{A-assmp: polynomial lipschitz smoothness} on the drift coefficient, in combination with Assumption \ref{assmp5: monotonicity}, are found in much of the tamed/ truncated scheme literature, for instance \cite{articletamed, articledif, doi:10.1137/S0036142901389530}, and are naturally required when applying Langevin-based sampling methods to light tailed densities \cite{BROSSE20193638}. The A-Assumption \ref{A-assmp: coercivity} is common in the sampling literature \cite{articleula, pmlr-v65-raginsky17a}, and is implied in particular by Assumption \ref{assmp5: monotonicity} when $\mu<0$. In most examples of interest one would expect A-Assumption \ref{A-assmp: coercivity} to hold, for example the SDE \eqref{eq: example SDE superlinear} considered in the simulations section. \begin{prop}
Suppose that A-Assumptions \ref{A-assmp: polynomial lipschitz}, \ref{A-assmp: polynomial lipschitz smoothness} and \ref{A-assmp: coercivity} hold, as well as Assumption \ref{assmp5: monotonicity} for the main SDE \eqref{eq: main equation}. Then the conclusions of Theorem \ref{thm: main result}, Corollary \ref{cor: strong convergence} and Corollary \ref{cor: weak convergence} hold.
\end{prop}
\begin{proof}
 We choose $V(x):=\lvert x \rvert^\xi+1$, so that $V\geq 1$ and since $\xi\geq 2$ one has $V\in C^2(\mathbb{R}^d)$ as required. It suffices to show Assumptions 1-4 from Section \ref{subsec: assumptions} hold, since we have assumed Assumption \ref{assmp5: monotonicity}. Therefore 
    \begin{itemize}
        \item   One has that Assumption \ref{assmp1: Lipschitz like behaviour} follows trivially from A-Assumption \ref{A-assmp: polynomial lipschitz}.
        \item Likewise Assumption \ref{assmp2: higher smoothness} follows from A-Assumption \ref{A-assmp: polynomial lipschitz smoothness} and a Taylor expansion argument, see Lemma 4 in \cite{SABANIS201984} for example.
        \item For Assumption \ref{assmp3: moment assmp}, note that for any $\alpha, \beta, r>0$ and any $V$ one has
    \begin{align}\label{eq: generator calculation}
    &\mathcal{L}e^{\beta V^\alpha} +re^{\beta V^\alpha}\nonumber \\
    &=\langle \nabla (e^{\beta V^\alpha}), b\rangle +\frac{1}{2}\Delta (e^{\beta V^\alpha})+re^{\beta V^\alpha}\nonumber \\
    &=\alpha \beta e^{\beta V^\alpha}(\alpha \beta V^{\alpha-1}\langle \nabla V, b\rangle + V^{\alpha-1}\Delta V+ (\alpha-1)V^{\alpha-2}\lvert \nabla V \rvert^2\nonumber \\
    &+\alpha \beta V^{2(\alpha-1)}\lvert \nabla V \rvert^2+r).
\end{align}
Completing this calculation will end up yielding a polynomial in $\lvert x \rvert$. Noting
\begin{equation}
  \nabla V=\xi \lvert x \rvert^{\xi-2}x,\;\;\;  \nabla V\leq c \lvert x \rvert^{\xi-2},
\end{equation}
and using A-Assumption \ref{A-assmp: coercivity}, one sees that the first term in the brackets of \eqref{eq: generator calculation} has leading order term $-c_1\lvert x \rvert^{\alpha \xi}$, where $c_1$ is a generic constant. All other terms have a positive sign and are of order $2\alpha \xi-2$ or lower. Therefore, bounding lower order terms by Young's inequality (i.e. $\lvert x \rvert^l\leq c(1+\lvert x \rvert^{l'})$ for $0<l<l'$) one may write for constants $c_1,c_2>0$ depending on $r>0$ that
\begin{align}\
    &\mathcal{L}e^{\beta V^\alpha} +re^{\beta V^\alpha}\leq  e^{(1+\lvert x \rvert^\xi)^\alpha}(c_1+c_1\lvert x \rvert^{2\alpha \xi-2}-c_2\lvert x \rvert^{\alpha \xi}).
\end{align}
Since the expression insides the brackets must be negative for sufficiently large $\lvert x \rvert$, the expression above must be bounded for every $r>0$ (not uniformly in $r>0$ however). Therefore Assumption \ref{assmp3: moment assmp} holds.
\item Finally Assumption \ref{assmp4: control of nabla V and norm} holds trivially for $w=1$.
    \end{itemize}

\end{proof}
\subsection{Lennard-Jones Type Coefficients}\label{subsec: lennard jones}
 In this section we consider coefficients that explode to infinity on a compact set, in the manner of the famous Lennard-Jones potential in molecular dynamics \cite{29acd3d494044594aea0829ef236aad6}. For example one may consider for $d=2$
\begin{equation}\label{eq: example SDE singular}
    dX_{1,t}=-\lvert X_{1,t}-X_{2,t}\rvert^{-\alpha-1}dt+dW_{1,t},\;\; dX_{2,t}=\lvert X_{1,t}-X_{2,t}\rvert^{-\alpha-1}dt+dW_{2,t}
\end{equation}
 for any $\alpha>1$ and $X^1_0<X^2_0$. For these dynamics, and dynamics of this kind, it's easy to show that for the Euler scheme $E[\lvert X^n_t\rvert ]=\infty$ for every $t>n^{-1}$, since for any one dimensional normal random variable of positive variance $E[\lvert \xi\rvert^{-\alpha}]=\infty$. 
 
 Specifically, we assume the following
 \begin{Bassumption}\label{B assmp}
     Let $ D:=\{x\in \mathbb{R}^d\vert x_1<x_2<...<x_d\}\subset \mathbb{R}^d$, and setting $\alpha>b>1$ and for $1\leq i<j\leq d$ coefficients $l_1,l_2,\nu\geq 0$, define
     \begin{equation}
         U(x)=\sum_{i=2}^{d}(l_1\lvert x_i-x_{i-1}\rvert^{-a}-l_2\lvert x_i-x_{i-1}\rvert^{-b})+\nu \lvert x\rvert^2
     \end{equation}
     and let the drift coefficient be given by
     \begin{equation}
         b:=-\nabla U.
     \end{equation}
 \end{Bassumption}
 Note we include only nearest-neighbour interactions for simplicity, more general potentials $U$ could be considered at the cost of more complicated calculations.
 \begin{prop}
Suppose that B-Assumption \ref{B assmp} and Assumption \ref{assmp5: monotonicity} hold for the main SDE \eqref{eq: main equation}. Then one has that the conclusions of Theorem \ref{thm: main result}, Corollary \ref{cor: strong convergence} and Corollary \ref{cor: weak convergence} hold.
\end{prop}
\begin{proof}
    We choose
    \begin{equation}
         V(x)=1+\sum_{i=2}^d\lvert x_i-x_{i-1}\rvert^{-a-3}+\frac{\nu}{2}\lvert x \rvert^2.
     \end{equation}
     and so that $V\geq 1$ and $V\in C^2(D)$ as required. Again we must verify Assumptions 1-4 from Section \ref{subsec: assumptions} hold.
     \begin{itemize}
         \item First let us note that $D\subset\mathbb{R}^d$ is convex, in that for $x,y\in D$, for any $t\in [0,1]$ one has $tx+(1-t)y\in D$. Therefore, we may use the higher dimensional fundamental theorem of calculus to yield
\begin{equation}
   \lvert b(x)-b(y)\rvert \leq \lvert x-y\rvert \sup_{t\in [0,1]}\lvert \nabla b(tx+(1-t)y\rvert ,
\end{equation}
and denoting the Frobenius norm of tensors on $\mathbb{R}^{d\times d\times d}$ as $\lvert \cdot\rvert$, one has furthermore
\begin{equation}
   \lvert b(x)-b(y)-\nabla b(x)(x-y)\rvert \leq \lvert x-y\rvert^2\sup_{t\in [0,1]}\lvert \nabla^2 b(tx+(1-t)y\rvert.
\end{equation}
Since $b=-\nabla U$, one may calculate
\begin{equation}
       \lvert   b(x)\rvert \leq c\sum_{i=2}^{d}(\lvert x_i-x_{i-1}\rvert^{-a-1}+\lvert x_i-x_{i-1}\rvert^{-b-1})+ c\lvert x \rvert,
     \end{equation}
     \begin{equation}
       \lvert   \nabla b(x)\rvert \leq c\sum_{i=2}^{d}(\lvert x_i-x_{i-1}\rvert^{-a-2}+\lvert x_i-x_{i-1}\rvert^{-b-2})+ c,
     \end{equation}
    \begin{equation}
       \lvert   \nabla^2 b(x)\rvert \leq c\sum_{i=2}^{d}(\lvert x_i-x_{i-1}\rvert^{-a-3}+\lvert x_i-x_{i-1}\rvert^{-b-3}),
     \end{equation} 
     Furthermore, for every $\gamma>0$, $i=2,...,d$ and $t\in [0,1]$ one has 
     \begin{equation}
     \lvert t(x_i-x_{i-1})+(1-t)(y_i-y_{i-1})\rvert^{-\gamma}\leq \lvert x_i-x_{i-1}\rvert^{-\gamma}+\lvert y_i-y_{i-1}\rvert^{-\gamma}.
     \end{equation}
     Combining these observations, and using Young's inequality to bound 
     \begin{equation}
     \lvert x_i-x_{i-1}\rvert^{-\gamma}\leq c(1+\lvert x_i-x_{i-1}\rvert^{-a-3}),
     \end{equation}
     for any $0<\gamma<a+3$, one sees that Assumptions \ref{assmp1: Lipschitz like behaviour} and \ref{assmp2: higher smoothness} follow.
     \item To show Assumption \ref{assmp3: moment assmp} let us begin by defining $\tilde{U}(x):=\sum_{i=2}^{d}\lvert x_i-x_{i-1}\rvert^{-a}$ and $\tilde{V}(x):=\sum_{i=2}^{d}\lvert x_i-x_{i-1}\rvert^{-a-3}$ and showing there exists a positive constant $c_1>0$ such that
     \begin{equation}\label{eq: coercivity bound}
         \langle \nabla \tilde{V}(x), \nabla \tilde{U}(x)\rangle \leq -c_1\sum_{i=2}^{d}\lvert x_i-x_{i-1}\rvert^{-a-5}.
     \end{equation}
     One notes by direct calculation that
\begin{align}\label{eq: V and U calculation}
   - \langle \nabla \tilde{V}(x), \nabla \tilde{U}(x)\rangle&=-2a (a+3)\sum_{i=2}^d\lvert x_i-x_{i-1}\rvert^{-2a-5}\nonumber \\
    &+a (a+3)\sum_{i=2}^{d-1}\lvert x_{i+1}-x_i\rvert^{-a-4}\lvert x_i-x_{i-1}\rvert^{-a-1}\nonumber \\
    &+a (a+3)\sum_{i=2}^{d-1}\lvert x_{i+1}-x_i\rvert^{-a-1}\lvert x_i-x_{i-1}\rvert^{-a-4}.
\end{align}
To complete the calculation we bound the latter two terms using the following result: for every $\gamma_1,\gamma_2>0$, $k\in \mathbb{N}$ there exists $\epsilon>0$ such that for every sequence $z_1,...,z_k>0$ one has
\begin{equation}\label{eq: sequence result}
    \sum_{i=1}^{k-1}z_{i+1}^{\gamma_1}z_i^{\gamma_2}\leq (1-\epsilon)\sum_{i=1}^kz_i^{\gamma_1+\gamma_2}.
\end{equation}
By a rescaling argument that it is sufficient to show 
\begin{equation}\label{eq: ratio bound}
\sup_{\sum_{i=1}^kz_i^{\gamma_1+\gamma_2}=1} \sum_{i=1}^{k-1}z_{i+1}^{\gamma_1}z_i^{\gamma_2}\leq 1-\epsilon.
\end{equation}
Now let $\sum_{i=1}^kz_i^{\gamma_1+\gamma_2}=1$ and let us bound by Young's inequality
\begin{align}
    \sum_{i=1}^{k-1}z_{i+1}^{\gamma_1}z_i^{\gamma_2}&\leq \frac{\gamma_1}{\gamma_1+\gamma_2}z_1^{\gamma_1+\gamma_2}+\frac{\gamma_2}{\gamma_1+\gamma_2}z_k^{\gamma_1+\gamma_2}+\sum_{i=2}^{k-1}z_i^{\gamma_1+\gamma_2}\nonumber \\
    &\leq \sum_{i=1}^kz_i^{\gamma_1+\gamma_2}.
\end{align}
Equality is achieved in the first inequality if $z_1=z_2=...=z_k$, and in the second inequality if $z_1=z_k=0$. Since these conditions are mutually exclusive under $\sum_{i=1}^kz_i^{\gamma_1+\gamma_2}=1$, for $z_1,...,z_k>0$ obeying $\sum_{i=1}^kz_i^{\gamma_1+\gamma_2}=1$ one must have
\begin{equation}\label{eq: ratio bound weak}
 \sum_{i=1}^{k-1}z_{i+1}^{\gamma_1}z_i^{\gamma_2}<1.
\end{equation}
One achieves \eqref{eq: ratio bound} from \eqref{eq: ratio bound weak} by compactness and continuity. Therefore, we have shown \eqref{eq: sequence result}, which we may apply to \eqref{eq: V and U calculation} to achieve \eqref{eq: coercivity bound}. From this, one sees by bounding lower order terms using Young's inequality, that there exists $c_1,c_2>0$ such that
\begin{align}
         - \langle \nabla V(x), \nabla U(x)\rangle &\leq -c_1\sum_{i=2}^{d}\lvert x_i-x_{i-1}\rvert^{-a-5}-2\nu^2\lvert x \rvert^2\nonumber \\
         &+c_2\sum_{i=2}^{d}\lvert x_i-x_{i-1}\rvert^{-a-4}+c_2.
     \end{align}
     Furthermore, one has by the triangle inequality
     \begin{equation}
         \lvert \nabla V(x)\rvert^2 \leq c\sum_{i=2}^{d}\lvert x_i-x_{i-1}\rvert^{-2a-8}+\nu^2\lvert x \rvert,
     \end{equation}
     \begin{equation}
         \lvert \Delta V(x)\rvert \leq c\sum_{i=2}^{d}\lvert x_i-x_{i-1}\rvert^{-a-5}.
     \end{equation}
     Now note by Jensen's inequality and Hölder's inequality, for every $\gamma>0$ there exists $c_1,c_2>0$ such that
     \begin{equation}\label{eq: power identity}
c_1V(x)^{\gamma/(a+3)}\leq \sum_{i=2}^{d}\lvert x_i-x_{i-1}\rvert^{-\gamma}\leq c_2V(x)^{\gamma/(a+3)}.
     \end{equation}
     Therefore, repeating the calculation from \eqref{eq: generator calculation}, one has for arbitrary $r>0$ and $\alpha>0$, and bounding lower order terms with Young's inequality, one has
     \begin{align}
      & \mathcal{L}e^{V^a} +re^{V^a}  =( V^{\alpha-1}\langle \nabla V, b\rangle + V^{\alpha-1}\Delta V+ (\alpha-1)V^{\alpha-2}\lvert \nabla V \rvert^2\nonumber \\
    &+V^{2(\alpha-1)}\lvert \nabla V \rvert^2+r)e^{\beta V^\alpha}\nonumber \\
    &\leq(-c_1 V^{\alpha+\frac{2}{\alpha+3}}+c_2 V^{\alpha+\frac{1}{\alpha+3}}+c_2-\nu^2\lvert x \rvert^2+r)e^{\beta V^\alpha}.
     \end{align}
     One sees this expression bounded for every $r>0$, at which point Assumption \ref{assmp3: moment assmp} holds. 
     \item Finally Assumption \ref{assmp4: control of nabla V and norm} holds by \eqref{eq: power identity}.
     \end{itemize}
\end{proof}

\section{Simulations}\label{sec: simulations}
In this section we present simulations on the rate of explosion, i.e. the frequency of trajectories for which the Euler scheme \eqref{eq: Euler scheme} satisfies $\sup_{t\in [0,1]}\lvert X^n_{\kappa_n(t)}\rvert\geq Cn^a$ for some chosen $C,a>0$. We look into the two forms of singular behaviour mentioned above: superlinearities that grow faster than linearly at infinity, and Lennard-Jones type singular behaviour where the coeffifients explode to infinity on a compact set. For the first class of coefficients we take as a representative example $b(x)=-x^3$ as in \eqref{eq: example SDE superlinear}, for the second we consider $b(x,y)=(-\lvert x-y\rvert^{-13},\lvert x-y\rvert^{-13})$ and initial conditions $X_{1,0}<X_{2,0}$ as in \eqref{eq: example SDE singular}.

Our experiments show that in the case of superlinear coefficients, there is a threshhold over which the probability of an explosion becomes extremely unlikely, and in the $10^5$ Monte-Carlo simulations we perform not a single explosion is seen. The value of $a$ in the threshhold $Cn^a$ seems not to matter here. On the other hand, for singular coefficients of the kind $b(x,y)=(-\lvert x-y\rvert^{-13},\lvert x-y\rvert^{-13})$, blowup behaviour is more subtle, and blowups were observed for even small stepsizes. Here $a>0$ does appear to matter, and the exponential decay can be observed in the graphs below. This perhaps motivates the very small stepsizes taken in particle simulations.

We conclude by noting that the problem of choosing $C,a>0$ is not obvious. The aim of this paper is to demonstrate that one can remove outliers, and virtually any method for doing this is asymptotically valid. However, asymptotic guarantees do not correspond to accurate performance in practise, and taming method can lead to more robust performance when blow-up frequency is not known. It should also be observed that in cases where tamed coefficients often differ from the true coefficients, i.e. where large trajectories are common, the approximation on these trajectories is unlikely to be accurate. Comparing taming method performance on exlposive trajectories is an interesting topic for future work.

\subsection{Superlinear Example}
We consider approximating SDE \eqref{eq: main equation} with $b(x)=-x^3$, i.e. \eqref{eq: example SDE superlinear}, with $X_0=6$. We used $10^5$ approximations for each value of $n$ and $a>0$ to produce estimates for the blow-up rate of the Euler scheme \eqref{eq: Euler scheme}, that is $P(\sup_{t\in [0,1]}\lvert X^n_{\kappa_n(t)}\rvert\geq 10 n^a)$. We consider parameter values $a=1/2,1,2$ and $n=10,...,20$. We note that for all $a=1/2,1,2$ and $n=10,...,16$ a blow-up was observed on each trajectory, whilst for $a=1/2,1,2$ and $n=20$ not a single blow-up was observed.

\begin{table}[H]
\caption{Estimates of $P(\sup_{t\in [0,1]}\lvert X^n_{\kappa_n(t)}\rvert\geq 10n^a)$ for $b(x)=-x^3$}
\centering
\begin{tabular}{l|llll}
n  &  & $a=1/2$ & $a=1$ & $a=2$ \\ \hline
10 &  &  1.0000 & 1.0000 & 1.0000\\
11 &  &  1.0000 & 1.0000 & 1.0000\\
12 &  & 1.0000 & 1.0000 & 1.0000 \\
13&  &  1.0000 & 1.0000 & 1.0000 \\
14 &  &  1.0000 & 1.0000 & 1.0000\\
15 &  &  1.0000 & 1.0000 & 1.0000 \\
16 &  &  1.0000 & 1.0000 & 1.0000\\
17 &  &  0.9998 & 0.9998 & 0.9998\\
18 &  &  0.5020 & 0.5011 & 0.5053\\
19 &  & 0.0003 & 0.0004 & 0.0004 \\
20&  &  0.0000 & 0.0000 & 0.0000\\
\end{tabular}
\end{table}
\begin{figure}[H]
\includegraphics[width=11cm]{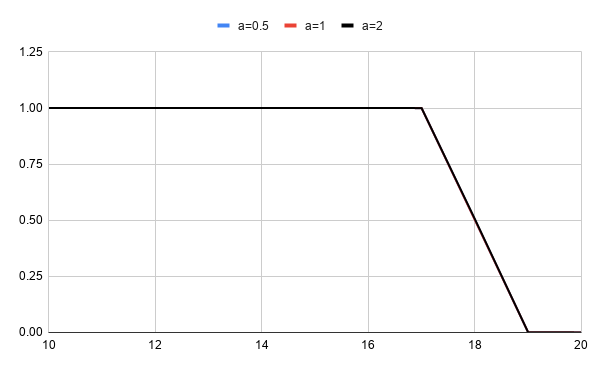}
\centering
\end{figure}

\subsection{Lennard-Jones Type Coefficients}
Here we consider $b(x,y)=(-\lvert x-y\rvert^{-13},\lvert x-y\rvert^{-13})$ with $X_0=(-1,1)$. We estimate $P(\sup_{t\in [0,1]}\lvert X^n_{\kappa_n(t)}\rvert\geq 1000 n^a)$ for $n=2,...,32$ and $a=1/10,1,2$. Here one may observe a much smoother exponential curve, suggesting the problem of choosing a sufficiently small stepsize is more subtle for coefficients that explode to infinity on a compact set.
\begin{table}[H]
\caption{Estimates of $P(\sup_{t\in [0,1]}\lvert X^n_{\kappa_n(t)}\rvert\geq 1000n^a)$ for $b(x,y)=(-\lvert x-y\rvert^{-13},\lvert x-y\rvert^{-13})$}
\centering
\begin{tabular}{l|llll}
n  &  & $a=1/10$ & $a=1$ & $a=2$ \\ \hline
2 &  & 0.07043 & 0.0662 & 0.06437\\
3 &  & 0.08222 & 0.07599 & 0.06727 \\
4&  &  0.08272 & 0.07004 & 0.06085 \\
5 &  &  0.0771 & 0.06468 & 0.05332\\
6 &  &  0.07192 & 0.05734 & 0.04566 \\
7 &  &  0.06616 & 0.05069 & 0.03933\\
8 &  &  0.05981 & 0.03982 & 0.03272\\
9 &  &  0.05477 & 0.03982 & 0.02792\\
10 &  & 0.0499 & 0.03389 & 0.02375 \\
11 &  &  0.04731 & 0.02986 & 0.02068\\
12 &  & 0.04109 & 0.02638 & 0.0166 \\
13&  &  0.03687 & 0.0228 & 0.01436 \\
14 &  &  0.03473 & 0.01982 & 0.01265\\
15 &  &  0.03114 & 0.01834 & 0.00979 \\
16 &  &  0.02866 & 0.01603 & 0.00831\\
17 &  &  0.02619 & 0.01327 & 0.00736\\
18 &  &  0.02369 & 0.01214 & 0.00665\\
19 &  & 0.02172 & 0.01042 & 0.00527 \\
20 &  & 0.01926 & 0.00908 & 0.00425 \\
21 &  &  0.0184 & 0.00813 & 0.00361\\
22 &  & 0.01686 & 0.00716 & 0.00289 \\
23&  &  0.01566 & 0.0068 & 0.00295 \\
24 &  &  0.01449 & 0.00537& 0.00227\\
25 &  & 0.0128 & 0.00516 & 0.00196 \\
26 &  &  0.0111 & 0.0042 & 0.00186\\
27 &  &  0.01044 & 0.00369 & 0.00132\\
28 &  &  0.00957 & 0.00355 & 0.0014\\
29 &  & 0.00919 & 0.00296 & 0.00101 \\
30 &  &  0.00875 & 0.00274 & 0.00098\\
31 &  &  0.00764 & 0.00235& 0.00068\\
32 &  & 0.00728 & 0.00209 & 0.00077 \\
\end{tabular}
\end{table}

\begin{figure}[H]
\includegraphics[width=11cm]{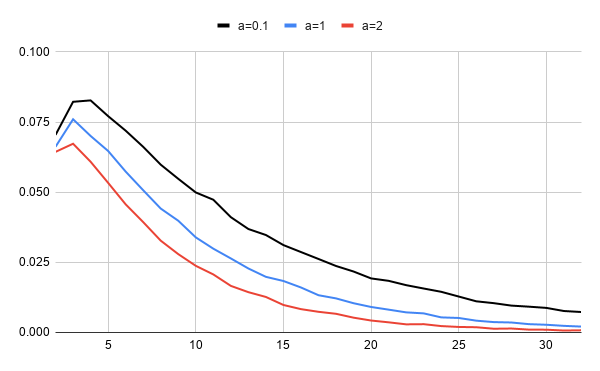}
\centering
\end{figure}



 \bibliographystyle{plain}

\bibliography{ref}

\end{document}